\documentclass[11pt,a4paper]{amsart}
\usepackage{amsmath,amssymb,amsthm,mathtools}
\usepackage[a4paper,margin=28mm]{geometry}
\usepackage{enumitem}
\usepackage[unicode,hidelinks]{hyperref}
\numberwithin{equation}{section}
\allowdisplaybreaks[2]
\setlist[enumerate]{leftmargin=2.2em,itemsep=2pt,topsep=4pt}
\newcommand{\T}{\mathbb T}
\newcommand{\R}{\mathbb R}
\newcommand{\C}{\mathbb C}
\newcommand{\Z}{\mathbb Z}
\newcommand{\dd}{\,\mathrm d}
\newcommand{\norm}[1]{\left\lVert #1\right\rVert}
\newcommand{\ip}[2]{\left\langle #1,#2\right\rangle}
\newcommand{\cH}{\mathcal H}
\DeclareMathOperator{\supp}{supp}

\usepackage[T1]{fontenc}
\usepackage[utf8]{inputenc}
\usepackage{lmodern}
\usepackage{microtype}
\theoremstyle{plain}
\newtheorem{theorem}{Theorem}[section]
\newtheorem{lemma}[theorem]{Lemma}
\newtheorem{proposition}[theorem]{Proposition}
\title[Endpoint exclusion for bounded damping]{Endpoint Exclusion for Bounded Damping on the Standard Flat Two-Torus}
\author{Zikang Deng}
\address{Beijing Normal University}
\date{}
\hypersetup{pdftitle={Endpoint Exclusion for Bounded Damping on the Standard Flat Two-Torus},pdfauthor={Zikang Deng}}
\begin{document}
\begin{abstract}
We prove that, for every fixed nonzero nonnegative bounded damping $W$ on the standard flat two-torus, the stationary damped wave operator satisfies $\|(-\Delta-s^2+isW)^{-1}\|_{L^2\to L^2}=o(|s|)$. Consequently, for initial data in $H^2\times H^1$, the uniform decay rate of the square root of the energy improves from $O(t^{-1/2})$ to $o(t^{-1/2})$. This excludes the endpoint quasimodes proposed by Kleinhenz in 2024. The proof combines a uniform $L^4$ estimate, equidistribution of densities when the largest Fourier coefficient tends to zero, and a one-dimensional equation for a nonzero modulated weak limit. No continuity of the damping is required. The improvement may depend on $W$ and does not assert a fixed gain in the polynomial exponent.
\end{abstract}
\maketitle

\section{Introduction and main results}\label{sec:intro}
\subsection{Background}
Let $\T=\R/(2\pi\Z)$ and equip $X=\T^2$ with the standard flat metric and normalized measure $\dd m=(2\pi)^{-2}\dd x_1\dd x_2$. We consider the damped wave equation
\begin{equation}\label{eq:wave}
 \partial_t^2u-\Delta u+W(x)\partial_tu=0,
 \qquad \Delta=\partial_{x_1}^2+\partial_{x_2}^2,
\end{equation}
where throughout the paper
\begin{equation}\label{eq:dampingassumptions}
 0\le W\in L^\infty(X),\qquad W\not\equiv0.
\end{equation}
Nonvanishing is understood modulo null sets, so that $\int_XW\dd m>0$. We assume neither continuity of $W$ nor that $\{W>0\}$ contains an open set. The associated energy is
\begin{equation}\label{eq:energy}
 E(u,t)=\frac12\int_X\bigl(|\nabla u(t,x)|^2+|\partial_tu(t,x)|^2\bigr)\dd m(x).
\end{equation}
For sufficiently regular solutions, integration by parts gives
\begin{equation}\label{eq:energydiss}
 \frac{\dd}{\dd t}E(u,t)=-\int_XW|\partial_tu|^2\dd m\le0.
\end{equation}
Anantharaman and L\'eautaud proved that Schr\"odinger observability implies $O(t^{-1/2})$ decay of the square root of the energy, and asked whether this rate is optimal in abstract systems and torus models \cite[Section 2C, Question (3)]{AL}. Together with the rough observability result of Burq and Zworski \cite{BZ}, this yields the same rate for every damping satisfying \eqref{eq:dampingassumptions} on the two-torus. We prove a strict little-$o$ improvement for each fixed damping on the standard flat two-torus.

Define the stationary operator
\begin{equation}\label{eq:Ps}
 P_W(s)=-\Delta-s^2+isW:H^2(X)\longrightarrow L^2(X),\qquad s\in\R.
\end{equation}
In his CIRM 2024 problem \cite{CIRM}, Kleinhenz asked for a fixed damping $W$, a constant $C>0$, frequencies $s_j\to\infty$, and functions $v_j\in H^2(X)$ such that
\begin{equation}\label{eq:CIRM}
 \norm{v_j}_2=1,\qquad \norm{P_W(s_j)v_j}_2\le C/s_j.
\end{equation}
We prove that such an endpoint quasimode sequence does not exist, and convert the resulting resolvent improvement into a uniform time-decay estimate.

\subsection{Notation and the main theorem}
All $L^p$ norms use normalized measure, and $\ip fg=\int_Xf\overline g\dd m$ is linear in its first argument. For $k\in\Z^2$, write $e_k(x)=e^{ik\cdot x}$ and $\widehat f(k)=\ip f{e_k}$. Thus
\begin{equation}\label{eq:parseval}
 f=\sum_{k\in\Z^2}\widehat f(k)e_k\ \text{in }L^2,
 \qquad \norm f_2^2=\sum_{k\in\Z^2}|\widehat f(k)|^2,
 \qquad -\Delta e_k=|k|^2e_k.
\end{equation}
We use the Sobolev norms $\norm f_{H^r}^2=\sum_k(1+|k|^2)^r|\widehat f(k)|^2$. Write $\dot H^r=H^r(X)/\C$ for $r=1,2$, and let $[u]$ denote the displacement modulo constants. We equip $\dot H^1$ with the norm $\norm{\nabla u}_2$ and $\dot H^2$ with the $H^2$ norm of the mean-zero representative. Set
\begin{equation}\label{eq:datanorm}
 \mathcal N(u_0,u_1)=\norm{[u_0]}_{\dot H^2}+\norm{u_1}_{H^1}.
\end{equation}
We write $\rightharpoonup$ for weak convergence and $\lfloor a\rfloor$ for the greatest integer not exceeding $a$. Unless indicated otherwise, constants may depend on the fixed damping $W$, but not on sequence indices or asymptotic parameters. Circle measures are also normalized to have total mass one.

\begin{theorem}[Endpoint exclusion]\label{thm:main}
Let $X=(\R/(2\pi\Z))^2$ and assume \eqref{eq:dampingassumptions}. Then:
\begin{enumerate}[label=(\roman*)]
\item For every real $s\ne0$, the operator $P_W(s):H^2(X)\to L^2(X)$ is invertible, and
\begin{equation}\label{eq:mainres}
 \lim_{|s|\to\infty}\frac{\norm{P_W(s)^{-1}}_{L^2\to L^2}}{|s|}=0.
\end{equation}
\item There is no quasimode sequence satisfying \eqref{eq:CIRM}.
\item There exists a function $\varepsilon_W:[1,\infty)\to[0,\infty)$, depending only on $W$, with $\varepsilon_W(t)\to0$, such that every solution of \eqref{eq:wave} with initial data $(u_0,u_1)\in H^2(X)\times H^1(X)$ satisfies
\begin{equation}\label{eq:mainenergy}
 E(u,t)^{1/2}\le\varepsilon_W(t)t^{-1/2}\mathcal N(u_0,u_1),\qquad t\ge1.
\end{equation}
\end{enumerate}
\end{theorem}
For fixed $W$, estimate \eqref{eq:mainenergy} is uniform in the initial data. No common improvement function is asserted for different dampings. A little-$o$ improvement also need not yield $O(t^{-1/2-\delta})$ for any $\delta>0$; the rate $t^{-1/2}/\log(2+t)$ illustrates this distinction.

\subsection{Outline of the proof}
The proof has seven stages. Section~\ref{sec:step1} reduces the resolvent estimate to the exclusion of quasimodes. Section~\ref{sec:step2} derives small dissipation and a bounded free residual from the quasimode equation. Section~\ref{sec:step3} localizes the Fourier mass to a shell of bounded squared-frequency width and establishes a uniform $L^4$ bound. Section~\ref{sec:step4} shows that a vanishing largest Fourier coefficient forces the densities to converge weakly to the constant one, contradicting small dissipation. Section~\ref{sec:step5} modulates around a nonvanishing Fourier coefficient and extracts a nonzero weak limit whose Fourier support has at most one direction. Section~\ref{sec:step6} excludes this limit using the averaged damping and a one-dimensional equation. Section~\ref{sec:step7} passes from the scalar resolvent to the wave generator and uses a semigroup decay theorem with a positive-increase majorant to retain the little-$o$ improvement.

\section{Reduction to quasimodes}\label{sec:step1}
\subsection{Invertibility at nonzero frequencies}
\begin{lemma}[Zeros of trigonometric polynomials]\label{lem:trigzeros}
A trigonometric polynomial on the one- or two-dimensional torus that vanishes on a set of positive measure is identically zero.
\end{lemma}
\begin{proof}
In one dimension, if $p(\theta)=\sum_{n=-N}^Na_ne^{in\theta}$, then
\[
 e^{iN\theta}p(\theta)=\sum_{n=-N}^Na_nz^{n+N},\qquad z=e^{i\theta}.
\]
The right-hand side is an ordinary polynomial and has only finitely many zeros unless it vanishes identically. In two dimensions, write
\[
 p(x_1,x_2)=\sum_{n=-N}^Nb_n(x_2)e^{inx_1}.
\]
If the zero set has positive area, Fubini's theorem gives a positive-measure set of $x_2$ for which the corresponding $x_1$ section has a positive-measure zero set. The one-dimensional assertion makes each such section identically zero. Every $b_n$ therefore vanishes on that set of $x_2$, and a second application of the one-dimensional assertion gives $b_n\equiv0$.
\end{proof}

\begin{proposition}\label{prop:invertibility}
For every real $s\ne0$, the operator $P_W(s):H^2(X)\to L^2(X)$ is a bounded isomorphism. As a closed operator on $L^2$ with domain $H^2$, it satisfies
\begin{equation}\label{eq:adjoint}
 P_W(s)^*=P_W(-s),\qquad
 \norm{P_W(s)^{-1}}_{2\to2}=\norm{P_W(-s)^{-1}}_{2\to2}.
\end{equation}
\end{proposition}
\begin{proof}
If $P_W(s)v=0$, integration by parts yields
\[
 0=\ip{P_W(s)v}v=\norm{\nabla v}_2^2-s^2\norm v_2^2+is\int_XW|v|^2\dd m.
\]
Taking imaginary parts and using $s\ne0$ and $W\ge0$ gives $Wv=0$ almost everywhere, hence $(-\Delta-s^2)v=0$. In Fourier variables,
\[
 (|k|^2-s^2)\widehat v(k)=0,\qquad k\in\Z^2.
\]
Only finitely many frequencies can occur, so $v$ is a trigonometric polynomial. It vanishes on the positive-measure set $\{W>0\}$, and Lemma~\ref{lem:trigzeros} gives $v=0$.

Set $A_0=-\Delta+1$ and $B_s=-s^2-1+isW$. The map $A_0:H^2\to L^2$ is invertible, and its inverse is the Fourier multiplier $(1+|k|^2)^{-1}$. Truncation to $|k|\le N$ has operator-norm error at most $(1+N^2)^{-1}$, so $A_0^{-1}$ is compact on $L^2$. Since $B_s$ is bounded,
\[
 P_W(s)A_0^{-1}=I+B_sA_0^{-1}
\]
is a compact perturbation of the identity. The Fredholm alternative upgrades injectivity to surjectivity, and the bounded inverse theorem gives the isomorphism. Finally, $-\Delta-s^2$ is self-adjoint and $(isW)^*=-isW$; the adjoint formula for bounded perturbations proves \eqref{eq:adjoint}.
\end{proof}

\subsection{The contradiction sequence}
\begin{proposition}\label{prop:quasimode}
As $s\to+\infty$, the estimate $\norm{P_W(s)^{-1}}=o(s)$ holds if and only if there is no sequence satisfying \eqref{eq:CIRM}.
\end{proposition}
\begin{proof}
If the little-$o$ estimate fails, there are $\varepsilon_0>0$ and $s_j\to\infty$ with $\norm{P_W(s_j)^{-1}}\ge\varepsilon_0s_j$. By the definition of the operator norm, choose
\[
 \norm{F_j}_2=1,\qquad z_j=P_W(s_j)^{-1}F_j,
 \qquad \norm{z_j}_2\ge\varepsilon_0s_j/2.
\]
Setting $v_j=z_j/\norm{z_j}_2$ and $f_j=F_j/\norm{z_j}_2$ gives
\[
 P_W(s_j)v_j=f_j,\qquad \norm{v_j}_2=1,
 \qquad \norm{f_j}_2\le2/(\varepsilon_0s_j).
\]
Conversely, \eqref{eq:CIRM} implies
\[
 1\le\norm{P_W(s_j)^{-1}}\norm{P_W(s_j)v_j}_2
 \le\norm{P_W(s_j)^{-1}}C/s_j,
\]
which rules out the little-$o$ estimate.
\end{proof}

Suppose henceforth that such a sequence exists. After discarding finitely many terms, write
\begin{equation}\label{eq:contradiction}
\begin{gathered}
 s_j\ge1,\quad s_j\to\infty,\quad \norm{v_j}_2=1,
 \quad f_j=P_W(s_j)v_j,\quad \norm{f_j}_2\le C_0/s_j,\\
 D_j=-\Delta-s_j^2,\qquad D_jv_j+is_jWv_j=f_j.
\end{gathered}
\end{equation}
Here both $W$ and $C_0$ are independent of $j$.

\section{Dissipation and the free residual}\label{sec:step2}
\begin{lemma}\label{lem:basic}
Let $g_j=s_jWv_j$. Then
\begin{align}
 \int_XW|v_j|^2\dd m&\le C_0s_j^{-2},\label{eq:small-damping}\\
 \norm{g_j}_2^2&\le C_0\norm W_\infty,\label{eq:g-bound}\\
 D_jv_j=f_j-ig_j,\qquad \norm{D_jv_j}_2&\le K,\label{eq:free-bound}
\end{align}
where $K=C_0+(C_0\norm W_\infty)^{1/2}$.
\end{lemma}
\begin{proof}
Taking the imaginary part of the inner product of \eqref{eq:contradiction} with $v_j$ gives
\[
 s_j\int_XW|v_j|^2\dd m
 =\operatorname{Im}\ip{f_j}{v_j}
 \le\norm{f_j}_2\norm{v_j}_2\le C_0/s_j,
\]
which proves \eqref{eq:small-damping}. Since $W^2\le\norm W_\infty W$,
\[
 \norm{g_j}_2^2=s_j^2\int_XW^2|v_j|^2\dd m
 \le s_j^2\norm W_\infty\int_XW|v_j|^2\dd m
 \le C_0\norm W_\infty.
\]
The identity $D_jv_j=f_j-ig_j$, the triangle inequality, and $s_j\ge1$ give \eqref{eq:free-bound}.
\end{proof}

In particular, for every fixed $a>0$,
\begin{equation}\label{eq:positive-area}
 a\int_{\{W\ge a\}}|v_j|^2\dd m\le C_0s_j^{-2}.
\end{equation}
Although $Wv_j\to0$ in $L^2$, the sequences $g_j=s_jWv_j$ and $D_jv_j$ are only known to be bounded. This bounded residual must be retained in the limiting equation.

\section{Frequency concentration and a uniform \texorpdfstring{$L^4$}{L4} estimate}\label{sec:step3}
\subsection{Squared-frequency shells}
Write $c_j(k)=\widehat v_j(k)$. Parseval's identity and \eqref{eq:free-bound} give
\begin{equation}\label{eq:shellbound}
 \sum_{k\in\Z^2}|c_j(k)|^2=1,\qquad
 \sum_{k\in\Z^2}(|k|^2-s_j^2)^2|c_j(k)|^2\le K^2.
\end{equation}
For $R>0$, define
\[
 S_{j,R}=\{k\in\Z^2:\bigl||k|^2-s_j^2\bigr|\le R\}.
\]
For $k\notin S_{j,R}$, the weight $(|k|^2-s_j^2)^2$ exceeds $R^2$. Hence
\[
 R^2\sum_{k\notin S_{j,R}}|c_j(k)|^2
 \le\sum_{k\notin S_{j,R}}(|k|^2-s_j^2)^2|c_j(k)|^2\le K^2,
\]
and thus
\begin{equation}\label{eq:shelltail}
 \sum_{k\notin S_{j,R}}|c_j(k)|^2\le K^2/R^2.
\end{equation}
Since $|k|^2$ is an integer, the interval $[s_j^2-R,s_j^2+R]$ contains at most $L_R:=\lfloor2R\rfloor+2$ distinct energy levels. This bounds the number of levels, not the number of lattice points on each circle.

\subsection{An \texorpdfstring{$L^2$}{L2} bound for the densities}
The single-level part of the following estimate is the classical $L^4$ estimate for circles in two dimensions; see \cite{Cooke,Zygmund}. We include the proof to make uniformity in the spectral center explicit.
\begin{lemma}\label{lem:l4bound}
There is an absolute constant $C_4$ such that, for every $\rho\in\R$ and $v\in H^2(X)$,
\begin{equation}\label{eq:l4bound}
 \norm v_4\le C_4\bigl(\norm v_2+\norm{(-\Delta-\rho)v}_2\bigr).
\end{equation}
\end{lemma}
\begin{proof}
First let $\phi=\sum_{|k|^2=m}a_ke_k$. Since $\norm\phi_4^4=\norm{\phi^2}_2^2$, consider
\[
 \widehat{\phi^2}(q)=\sum_{\substack{|k|^2=m\\|q-k|^2=m}}a_ka_{q-k}.
\]
If $q\ne0$, subtracting the two circle equations gives $2q\cdot k=|q|^2$. This line meets $|k|^2=m$ in at most two points, so
\[
 |\widehat{\phi^2}(q)|^2
 \le2\sum_{\substack{|k|^2=m\\|q-k|^2=m}}|a_k|^2|a_{q-k}|^2.
\]
For $q=0$, Cauchy--Schwarz gives
\[
 |\widehat{\phi^2}(0)|^2
 =\left|\sum_{|k|^2=m}a_ka_{-k}\right|^2
 \le\left(\sum_{|k|^2=m}|a_k|^2\right)
     \left(\sum_{|k|^2=m}|a_{-k}|^2\right)=\norm\phi_2^4.
\]
On summing over nonzero $q$, each ordered pair $(k,\ell)$ occurs only at $q=k+\ell$. Consequently,
\[
 \sum_{q\ne0}|\widehat{\phi^2}(q)|^2
 \le2\sum_{\substack{|k|^2=m\\|\ell|^2=m}}|a_k|^2|a_\ell|^2
 =2\norm\phi_2^4.
\]
Parseval's identity now yields
\begin{equation}\label{eq:singlel4}
 \norm\phi_4^4\le3\norm\phi_2^4.
\end{equation}

Decompose into nonnegative integer energy levels:
\[
 v=\sum_{m\ge0}v_m,\qquad v_m=\sum_{|k|^2=m}\widehat v(k)e_k.
\]
For a finite sum, \eqref{eq:singlel4} and weighted Cauchy--Schwarz imply
\begin{align*}
 \norm v_4&\le3^{1/4}\sum_{m\ge0}\norm{v_m}_2\\
 &\le3^{1/4}\left(\sum_{m\ge0}\frac1{1+(m-\rho)^2}\right)^{1/2}
 \left(\sum_{m\ge0}(1+(m-\rho)^2)\norm{v_m}_2^2\right)^{1/2}.
\end{align*}
Choose an integer $n_0$ with $|\rho-n_0|\le1/2$ and set $a=\rho-n_0$. Since $|n-a|\ge|n|/2$ for $n\ne0$,
\[
 \sum_{m\ge0}\frac1{1+(m-\rho)^2}
 \le\sum_{n\in\Z}\frac1{1+(n-a)^2}
 \le1+8\sum_{n=1}^\infty n^{-2}=:C_*.
\]
Orthogonality of the energy levels gives
\[
 \sum_{m\ge0}(1+(m-\rho)^2)\norm{v_m}_2^2
 =\norm v_2^2+\norm{(-\Delta-\rho)v}_2^2.
\]
This proves the estimate for finite sums. For general $v\in H^2$, the same weighted estimate gives $\sum_m\norm{v_m}_4<\infty$. The partial sums therefore converge in $L^4$; since their $L^2$ limit is $v$, their $L^4$ limit is also $v$. This proves \eqref{eq:l4bound}.
\end{proof}

Applying the lemma to the contradiction sequence with $\rho=s_j^2$ gives
\[
 \sup_j\norm{v_j}_4\le C_4(1+K),\qquad
 \sup_j\norm{|v_j|^2}_2\le C_4^2(1+K)^2.
\]
In particular, for every measurable $F\subset X$,
\[
 \int_F|v_j|^2\dd m\le C_4^2(1+K)^2m(F)^{1/2}.
\]
The same bound permits approximation of measurable test functions by trigonometric polynomials in $L^2$:
\[
 \left|\int_X(H-H_N)|v_j|^2\dd m\right|
 \le\norm{H-H_N}_2\norm{v_j}_4^2.
\]

\section{Equidistribution and a nonvanishing coefficient}\label{sec:step4}
\subsection{An equidistribution lemma}
\begin{lemma}\label{lem:uniformdensity}
Suppose $u_j\in H^2(X)$, $\rho_j\to\infty$, and
\[
 \norm{u_j}_2=1,\qquad \norm{(-\Delta-\rho_j)u_j}_2\le K.
\]
If $b_j:=\sup_{k\in\Z^2}|\widehat u_j(k)|\to0$, then $|u_j|^2\rightharpoonup1$ in $L^2(X)$. Equivalently, for every $F\in L^2(X)$,
\begin{equation}\label{eq:uniformdensity}
 \int_XF|u_j|^2\dd m\longrightarrow\int_XF\dd m.
\end{equation}
\end{lemma}
\begin{proof}
Set $c_j(k)=\widehat u_j(k)$. Lemma~\ref{lem:l4bound} bounds the densities $|u_j|^2$ uniformly in $L^2$. Their Fourier coefficients are
\[
 \widehat{|u_j|^2}(q)=\sum_{k\in\Z^2}c_j(k+q)\overline{c_j(k)}.
\]
This series is absolutely convergent, since Cauchy--Schwarz gives $\sum_k|c_j(k+q)c_j(k)|\le1$. The convolution formula follows first for finite Fourier sums and then by $L^2$ convergence and $L^1$ convergence of the products. The zero-frequency coefficient is always $\widehat{|u_j|^2}(0)=1$.

Fix $q\ne0$ and $R>1$, and define
\[
 S_{j,R}=\{k:\bigl||k|^2-\rho_j\bigr|\le R\},\qquad
 a_{j,R}(k)=\mathbf1_{S_{j,R}}(k)c_j(k).
\]
The estimate used in \eqref{eq:shelltail} gives
\[
 \norm{c_j-a_{j,R}}_{\ell^2}\le K/R,\qquad
 \norm{c_j}_{\ell^2}=1,\qquad \norm{a_{j,R}}_{\ell^2}\le1.
\]
Abbreviate $c=c_j$ and $a=a_{j,R}$. The identity
\begin{align*}
 &\sum_kc(k+q)\overline{c(k)}-\sum_ka(k+q)\overline{a(k)}\\
 &\quad=\sum_k(c-a)(k+q)\overline{c(k)}
       +\sum_ka(k+q)\overline{(c-a)(k)}
\end{align*}
and Cauchy--Schwarz imply
\[
 \left|\widehat{|u_j|^2}(q)-\sum_ka_{j,R}(k+q)\overline{a_{j,R}(k)}\right|
 \le\norm{c-a}_{\ell^2}\norm c_{\ell^2}+\norm a_{\ell^2}\norm{c-a}_{\ell^2}
 \le2K/R.
\]

A nonzero term in the truncated sum requires $k,k+q\in S_{j,R}$. The integers $m=|k|^2$ and $n=|k+q|^2$ each have at most $L_R=\lfloor2R\rfloor+2$ possible values. For fixed $(m,n)$,
\[
 |k|^2=m,\qquad 2k\cdot q=n-m-|q|^2.
\]
A circle and a nondegenerate line intersect in at most two points, so all pairs $(m,n)$ contribute at most $2L_R^2$ lattice points. Each term has absolute value at most $b_j^2$, yielding
\begin{equation}\label{eq:densitycoefficient}
 \left|\widehat{|u_j|^2}(q)\right|\le2K/R+2L_R^2b_j^2,\qquad q\ne0.
\end{equation}
First let $j\to\infty$ with $R$ fixed, and then let $R\to\infty$. Every nonzero Fourier coefficient tends to zero. Since the zero-frequency coefficient is one, \eqref{eq:uniformdensity} follows for trigonometric polynomials.

For general $F\in L^2$, choose trigonometric polynomials $F_N\to F$ in $L^2$ and write
\[
 \int_XF(|u_j|^2-1)\dd m
 =\int_XF_N(|u_j|^2-1)\dd m
 +\int_X(F-F_N)(|u_j|^2-1)\dd m.
\]
The absolute value of the second term is at most
\[
 \norm{F-F_N}_2\left(\sup_j\norm{u_j}_4^2+1\right).
\]
Taking $j\to\infty$ first and then $N\to\infty$ completes the proof.
\end{proof}

\subsection{A nonvanishing Fourier coefficient}
If a subsequence of the contradiction sequence satisfied $\sup_k|c_j(k)|\to0$, Lemma~\ref{lem:uniformdensity} and $W\in L^2(X)$ would give
\[
 \int_XW|v_j|^2\dd m\longrightarrow\int_XW\dd m>0,
\]
contradicting \eqref{eq:small-damping}. Thus $\liminf_j\sup_k|c_j(k)|>0$. Choose $\delta>0$ and, after discarding finitely many terms, select $k_j\in\Z^2$ for each $j$ such that
\begin{equation}\label{eq:atom}
 |c_j(k_j)|\ge\delta.
\end{equation}

\section{Modulated limits and one-dimensional Fourier support}\label{sec:step5}
\subsection{A nonzero weak limit}
Multiplying $v_j,f_j,g_j$ by the same unimodular constant $\overline{c_j(k_j)}/|c_j(k_j)|$ preserves all equations and estimates. We may therefore assume $c_j(k_j)\ge\delta$. By \eqref{eq:shellbound},
\[
 \eta_j:=|k_j|^2-s_j^2,\qquad |\eta_j|\le K/\delta.
\]
Thus $|k_j|\to\infty$, and after taking a subsequence we may assume $\eta_j\to\eta\in\R$. Define the modulated sequences
\begin{equation}\label{eq:modulation}
 w_j=e^{-ik_j\cdot x}v_j,\qquad
 h_j=e^{-ik_j\cdot x}g_j=s_jWw_j,\qquad
 r_j=e^{-ik_j\cdot x}f_j.
\end{equation}
These functions remain periodic, and
\[
 \widehat w_j(q)=c_j(k_j+q),\qquad
 \norm{w_j}_2=1,\quad \sup_j\norm{h_j}_2<\infty,\quad
 \norm{r_j}_2\le C_0/s_j\to0.
\]
Weak compactness in $L^2(X)\times L^2(X)$ gives a further subsequence with
\begin{equation}\label{eq:weaklimits}
 w_j\rightharpoonup w,\qquad h_j\rightharpoonup h
 \quad\text{in }L^2(X).
\end{equation}
The zero-frequency coefficient passes to the weak limit, so
\begin{equation}\label{eq:nonzero-limit}
 \widehat w(0)=\lim_jc_j(k_j)\ge\delta>0.
\end{equation}
In particular, $w\ne0$. Taking the original equation at frequency $k_j+q$ and expanding $|k_j+q|^2$ gives
\begin{equation}\label{eq:modulatedmode}
 (\eta_j+2k_j\cdot q+|q|^2)\widehat w_j(q)
 +i\widehat h_j(q)=\widehat r_j(q),\qquad q\in\Z^2.
\end{equation}

\subsection{Structure of the frequency subgroup}
For fixed $q\in\Z^2$, the integer sequence $k_j\cdot q$ either takes some value infinitely often or tends to infinity in absolute value. Enumerating $\Z^2$, extracting nested subsequences, and taking a diagonal subsequence ensures simultaneously that, for each fixed $q$, the sequence $k_j\cdot q$ is eventually constant or $|k_j\cdot q|\to\infty$. The diagonal subsequence is eventually contained in each selected subsequence and preserves all preceding convergence properties.

Define
\begin{equation}\label{eq:Lambda}
 \Lambda=\{q\in\Z^2:k_j\cdot q\ \text{is eventually constant}\}.
\end{equation}
The set is closed under addition and negation, hence is a subgroup. If $mq\in\Lambda$ for a nonzero integer $m$, then $m(k_j\cdot q)$ is eventually constant, so $q\in\Lambda$; thus the subgroup is saturated. It cannot contain two linearly independent vectors: otherwise $k_j$ would eventually solve a fixed invertible linear system and would be constant, contrary to $|k_j|\to\infty$.

If $\Lambda\ne\{0\}$, divide a nonzero $p\in\Lambda$ by the greatest common divisor of its coordinates to obtain a primitive vector $q_0=(a,b)\in\Lambda$. Every integer vector collinear with $q_0$ is an integer multiple of it: if $q=tq_0\in\Z^2$, choose $\alpha a+\beta b=1$; then $t=\alpha q_1+\beta q_2\in\Z$. Consequently,
\begin{equation}\label{eq:rankone}
 \Lambda=\{0\}\quad\text{or}\quad\Lambda=\Z q_0,
\end{equation}
where $q_0$ is nonzero and primitive in the second case. If $q\notin\Lambda$, equation \eqref{eq:modulatedmode} gives
\[
 \widehat w_j(q)=\frac{\widehat r_j(q)-i\widehat h_j(q)}
 {\eta_j+2k_j\cdot q+|q|^2}\longrightarrow0.
\]
The numerator is uniformly bounded while the denominator tends to infinity in absolute value. Weak convergence therefore yields
\begin{equation}\label{eq:limitsupport}
 \supp\widehat w\subset\Lambda.
\end{equation}

\subsection{Damping constraints on the limits}
Modulation preserves absolute values, so \eqref{eq:small-damping} gives
\[
 \norm{W^{1/2}w_j}_2^2=\int_XW|v_j|^2\dd m\longrightarrow0.
\]
Multiplication by $W^{1/2}$ is bounded and preserves weak convergence. Hence
\begin{equation}\label{eq:Wwzero}
 W^{1/2}w=0,\qquad Ww=0\quad\text{almost everywhere}.
\end{equation}
On the other hand, let $Z=\{W=0\}$. Since $h_j=s_jWw_j$, we have $\mathbf1_Zh_j=0$. Passing weakly to the limit gives
\begin{equation}\label{eq:hz}
 h=0\quad\text{almost everywhere on }Z.
\end{equation}
If $\Lambda=\{0\}$, then $w$ is a nonzero constant, and \eqref{eq:Wwzero} forces $W=0$ almost everywhere, a contradiction. It remains to treat $\Lambda=\Z q_0$. The conditions \eqref{eq:Wwzero} and \eqref{eq:hz} constrain different limits: the growing factor $s_j$ in $h_j=s_jWw_j$ prevents a direct deduction of $h=0$ from \eqref{eq:Wwzero}.

\section{The one-dimensional limit equation and endpoint exclusion}\label{sec:step6}
\subsection{Periodic coordinates and averaged damping}
By \eqref{eq:limitsupport}, there exists $V\in L^2(\T)$ such that
\begin{equation}\label{eq:V}
 w(x)=V(q_0\cdot x),\qquad \widehat V(n)=\widehat w(nq_0).
\end{equation}
Existence also follows directly from $\sum_n|\widehat w(nq_0)|^2\le\norm w_2^2$. Write $q_0=(a,b)$, choose integers $c,d$ with $ad-bc=1$, and set
\[
 \theta=ax_1+bx_2,\qquad \phi=cx_1+dx_2\pmod{2\pi}.
\]
The inverse is $x_1=d\theta-b\phi$, $x_2=-c\theta+a\phi$. Both maps have integer coefficients and determinant one, so they define a torus automorphism preserving normalized measure. In these coordinates, $w=V(\theta)$. Denote the coordinate expressions of $W,h$ by $\widetilde W,\widetilde h$, and define
\begin{equation}\label{eq:averages}
 A(\theta)=\int_\T\widetilde W(\theta,\phi)\dd m(\phi),\qquad
 G(\theta)=\int_\T\widetilde h(\theta,\phi)\dd m(\phi).
\end{equation}
Fubini's theorem and Cauchy--Schwarz yield
\begin{equation}\label{eq:averagefacts}
 \begin{gathered}
 0\le A\le\norm W_\infty,\qquad \int_\T A\dd m=\int_XW\dd m>0,\\
 G\in L^2(\T),\qquad \norm G_2\le\norm h_2.
 \end{gathered}
\end{equation}
Indeed, $|G(\theta)|^2\le\int_\T|\widetilde h(\theta,\phi)|^2\dd m(\phi)$. Moreover,
\begin{align*}
 \widehat G(n)&=\int_{\T^2}\widetilde h(\theta,\phi)e^{-in\theta}\dd m(\theta)\dd m(\phi)\\
 &=\int_Xh(x)e^{-in q_0\cdot x}\dd m(x)=\widehat h(nq_0).
\end{align*}

\subsection{The limit equation and its regularity}
Since $q_0\in\Lambda$, there is a fixed integer $\ell$ such that, after discarding finitely many terms, $k_j\cdot q_0=\ell$. Taking $q=nq_0$ in \eqref{eq:modulatedmode} gives
\[
 (\eta_j+2n\ell+n^2|q_0|^2)\widehat w_j(nq_0)
 +i\widehat h_j(nq_0)=\widehat r_j(nq_0).
\]
For fixed $n$, weak convergence and $r_j\to0$ allow passage to the limit:
\begin{equation}\label{eq:1dFourier}
 (|q_0|^2n^2+2\ell n+\eta)\widehat V(n)=-i\widehat G(n).
\end{equation}
Equivalently, in the sense of distributions,
\begin{equation}\label{eq:1dequation}
 LV+iG=0,\qquad L=-|q_0|^2\partial_\theta^2-2i\ell\partial_\theta+\eta.
\end{equation}
Since $|q_0|^2>0$, there are $N\ge1$ and $c_*>0$ such that, for $|n|\ge N$,
\[
 \bigl||q_0|^2n^2+2\ell n+\eta\bigr|\ge c_*(1+n^2).
\]
Consequently,
\[
 \sum_{|n|\ge N}(1+n^2)^2|\widehat V(n)|^2
 \le c_*^{-2}\sum_{|n|\ge N}|\widehat G(n)|^2
 \le c_*^{-2}\norm G_2^2.
\]
The finitely many low frequencies are controlled by $V\in L^2$. Hence
\begin{equation}\label{eq:VH2}
 V\in H^2(\T).
\end{equation}
Thus \eqref{eq:1dequation} is an almost-everywhere identity between $L^2$ functions.

\subsection{Zero sets and elimination of the residual}
\begin{lemma}\label{lem:zeroderivatives}
If $V\in H^2(\T;\C)$, then $V'=V''=0$ almost everywhere on $\{V=0\}$.
\end{lemma}
\begin{proof}
A real-valued $h\in H^1(\T)$ has an absolutely continuous representative and is differentiable almost everywhere. Fix a level set $F=\{h=c\}$. Its isolated points are at most countable, since each can be isolated within $F$ by an interval with rational endpoints. Thus almost every $\theta\in F$ is both a differentiability point and an accumulation point of $F$. Choosing $\theta_n\in F\setminus\{\theta\}$ with $\theta_n\to\theta$, we obtain
\[
 h'(\theta)=\lim_n\frac{h(\theta_n)-h(\theta)}{\theta_n-\theta}=0.
\]
Applying this assertion to the real and imaginary parts of $V$ gives $V'=0$ almost everywhere on $\{V=0\}$. Applying it again to the real and imaginary parts of $V'\in H^1$ gives $V''=0$ almost everywhere on $\{V'=0\}$. Since $\{V=0\}\subset\{V'=0\}$ modulo a null set, the conclusion follows.
\end{proof}

Equations \eqref{eq:Wwzero} and \eqref{eq:V}, integrated in $\phi$, give
\begin{equation}\label{eq:AV}
 A(\theta)|V(\theta)|^2=0\quad\text{almost everywhere}.
\end{equation}
Set $B=\{A=0\}$ and $E=\{A>0\}$. On $B$,
\[
 \int_{B\times\T}\widetilde W(\theta,\phi)\dd m(\theta)\dd m(\phi)
 =\int_BA\dd m=0.
\]
Nonnegativity implies $\widetilde W=0$ almost everywhere on $B\times\T$, and \eqref{eq:hz} then gives $\widetilde h=0$ there. Hence
\begin{equation}\label{eq:GB}
 G=0\quad\text{almost everywhere on }B.
\end{equation}
On $E$, equation \eqref{eq:AV} gives $V=0$ almost everywhere. By \eqref{eq:VH2} and Lemma~\ref{lem:zeroderivatives}, we have $V=V'=V''=0$ almost everywhere on $E$, hence $LV=0$. Equation \eqref{eq:1dequation} yields
\begin{equation}\label{eq:GE}
 G=0\quad\text{almost everywhere on }E.
\end{equation}
Since $B\cup E=\T$, the two identities imply
\begin{equation}\label{eq:homogeneous}
 G=0,\qquad LV=0\quad\text{on }\T.
\end{equation}

\subsection{Proof of the stationary conclusions}
Equations \eqref{eq:1dFourier} and \eqref{eq:homogeneous} imply
\[
 (|q_0|^2n^2+2\ell n+\eta)\widehat V(n)=0,\qquad n\in\Z.
\]
The quadratic polynomial has at most two integer zeros, so $V$ is a trigonometric polynomial. By \eqref{eq:averagefacts}, the set $E=\{A>0\}$ has positive measure, and \eqref{eq:AV} says that $V$ vanishes there almost everywhere. Lemma~\ref{lem:trigzeros} gives $V\equiv0$, contradicting
\[
 \widehat V(0)=\widehat w(0)\ge\delta>0.
\]
Thus no sequence satisfying \eqref{eq:contradiction} exists. Proposition~\ref{prop:quasimode} yields the $o(s)$ resolvent estimate at positive frequencies, and the adjoint relation in Proposition~\ref{prop:invertibility} gives the same conclusion at negative frequencies. This proves parts (i) and (ii) of Theorem~\ref{thm:main}.

\section{Resolvent bounds and energy decay}\label{sec:step7}
We now deduce the time-decay estimate from the scalar bound
\begin{equation}\label{sg:scalarinput}
 N_W(s):=\norm{P_W(s)^{-1}}_{L^2\to L^2}=o(|s|),\qquad |s|\to\infty.
\end{equation}
This proves part (iii) of Theorem~\ref{thm:main}. The essential points are to obtain an $o(s^2)$ bound for the generator resolvent and to preserve the little-$o$ improvement in the semigroup argument.

\subsection{The generator on the reduced energy space}
Setting $v=\partial_tu$, we obtain
\begin{equation}\label{sg:firstorder}
 \partial_tu=v,\qquad \partial_tv=\Delta u-Wv.
\end{equation}
The energy is insensitive to constants in the displacement. We therefore quotient only the displacement and retain the constant velocity component. Define
\begin{equation}\label{sg:space}
 \cH_W=\dot H^1\times L^2,\qquad
 \norm{([u],v)}_{\cH_W}^2=\norm{\nabla u}_2^2+\norm v_2^2,
\end{equation}
\begin{equation}\label{sg:generator}
 G_W([u],v)=([v],\Delta u-Wv),\qquad D(G_W)=\dot H^2\times H^1.
\end{equation}
Poincar\'e's inequality makes the displacement norms equivalent to the quotient norms. The state $U=([u],v)$ satisfies $U'=G_WU$ and $E(u,t)=\frac12\norm{U(t)}_{\cH_W}^2$.

\begin{lemma}\label{sg:maximal}
The operator $G_W$ is densely defined, closed, and maximally dissipative. It generates a contraction $C_0$ semigroup $T_W(t)$ on $\cH_W$, and $0\in\rho(G_W)$.
\end{lemma}
\begin{proof}
Density of smooth functions gives density of the domain. Integration by parts yields
\begin{align}
 \operatorname{Re}\ip{G_W([u],v)}{([u],v)}_{\cH_W}
 &=\operatorname{Re}\left(\int\nabla v\cdot\nabla\overline u\dd m+
 \int(\Delta u-Wv)\overline v\dd m\right)\notag\\
 &=-\int W|v|^2\dd m\le0.\label{sg:dissipative}
\end{align}
To prove closedness, suppose $([u_n],v_n)\to([u],v)$ and $G_W([u_n],v_n)\to([a],b)$ in $\cH_W$. The convergences $v_n\to v$ in $L^2$ and $[v_n]\to[a]$ in $\dot H^1$ imply $v\in H^1$ and $[v]=[a]$. The second component gives $\Delta u=b+Wv\in L^2$ in distributions. Elliptic regularity yields $[u]\in\dot H^2$ and $G_W([u],v)=([a],b)$.

Fix $\lambda>0$ and $([f],g)\in\cH_W$, with $f$ the mean-zero representative. In solving $(\lambda-G_W)([u],v)=([f],g)$, adjust the constant in $u$ so that $v=\lambda u-f$. Elimination gives
\begin{equation}\label{sg:lambda}
 (-\Delta+\lambda^2+\lambda W)u=g+(\lambda+W)f.
\end{equation}
The associated sesquilinear form on $H^1$ is
\[
 a(u,z)=\int\bigl(\nabla u\cdot\nabla\overline z+
 \lambda^2u\overline z+\lambda Wu\overline z\bigr)\dd m,
\]
It is bounded and satisfies $\operatorname{Re}a(u,u)\ge\min\{1,\lambda^2\}\norm u_{H^1}^2$. Lax--Milgram gives a unique weak solution $u\in H^1$. Since the right-hand side lies in $L^2$ and $W\in L^\infty$, we have $\Delta u\in L^2$, hence $u\in H^2$ and $v\in H^1$. Thus $\lambda-G_W$ is surjective, and the Lumer--Phillips theorem gives the contraction semigroup.

Finally, solve $G_W([u],v)=([f],g)$. The first component gives $v=f+c$ and the second gives $\Delta u=g+W(f+c)$. Solvability is equivalent to a zero mean on the right-hand side, so choose
\begin{equation}\label{sg:zeroformula}
 c=-\frac{\int g\dd m+\int Wf\dd m}{\int W\dd m}.
\end{equation}
The denominator is positive, and the Poisson equation then has a unique mean-zero solution $u\in H^2$. For homogeneous data, integration gives $c\int W=0$, hence $c=0$ and $[u]=0$; thus $G_W$ is also injective. Formula \eqref{sg:zeroformula}, Poincar\'e's inequality, and elliptic estimates give
\[
 \norm{[u]}_{\dot H^2}+\norm v_{H^1}
 \le C_W\bigl(\norm{\nabla f}_2+\norm g_2\bigr).
\]
The inverse is therefore bounded.
\end{proof}

\subsection{Comparison of scalar and generator resolvents}
\begin{proposition}\label{sg:comparison}
We have $i\R\subset\rho(G_W)$. If $q_W(s)=\norm{(is-G_W)^{-1}}_{\cH_W\to\cH_W}$, then, for $|s|\ge1$,
\begin{equation}\label{sg:comparebound}
 q_W(s)\le C_W\bigl(1+|s|N_W(s)\bigr).
\end{equation}
In particular,
\begin{equation}\label{sg:generatorlittleo}
 q_W(s)=o(s^2),\qquad |s|\to\infty.
\end{equation}
\end{proposition}
\begin{proof}
The zero frequency is covered by Lemma~\ref{sg:maximal}. For $s\ne0$, take $f$ to be the mean-zero representative of $[f]$. Choosing a representative of $u$, the equation $(is-G_W)([u],v)=([f],g)$ becomes
\[
 v=isu-f,\qquad P_W(s)u=g+(is+W)f.
\]
Proposition~\ref{prop:invertibility} gives the unique solution
\begin{equation}\label{sg:eliminate}
 u=P_W(s)^{-1}\bigl(g+(is+W)f\bigr),\qquad v=isu-f.
\end{equation}
Here $u\in H^2$ and $v\in H^1$, and the homogeneous equation has only the zero state as a solution. The closed graph theorem gives a bounded inverse, so $i\R\subset\rho(G_W)$.

Assume $|s|\ge1$ and abbreviate $N=N_W(s)$ and $K_s=1+|s|N$. If $z=P_W(s)^{-1}h$, the real part of the inner-product identity gives
\begin{equation}\label{sg:scalarenergy}
 \norm{\nabla z}_2^2=s^2\norm z_2^2+\operatorname{Re}\ip h z
 \le(s^2N^2+N)\norm h_2^2.
\end{equation}
Since $\sqrt N\le(1+N)/2\le K_s/2$,
\begin{equation}\label{sg:gradientinverse}
 \norm{\nabla P_W(s)^{-1}h}_2\le(|s|N+\sqrt N)\norm h_2
 \le C K_s\norm h_2.
\end{equation}
By \eqref{eq:adjoint}, this estimate also holds for $P_W(-s)^{-1}$, and
\begin{equation}\label{sg:adjointinverse}
 (P_W(s)^{-1})^*=P_W(-s)^{-1},\qquad N_W(-s)=N_W(s).
\end{equation}

To avoid losing an extra frequency factor in a direct estimate of \eqref{sg:eliminate}, we use $f\in H^1$ and compute the velocity equation:
\begin{align}
 P_W(s)v&=isP_W(s)u-P_W(s)f\notag\\
 &=isg+is(is+W)f-(-\Delta f-s^2f+isWf)\notag\\
 &=\Delta f+isg.\label{sg:velocityeq}
\end{align}
Here $\Delta f$ is interpreted distributionally. For arbitrary $h\in L^2$, set $z=P_W(-s)^{-1}h\in H^2$. The adjoint identity and integration by parts give
\begin{align}
 \ip v h&=\ip v{P_W(-s)z}=\ip{P_W(s)v}z\notag\\
 &=\ip{\Delta f}z+is\ip g z
 =-\ip{\nabla f}{\nabla z}+is\ip g z.\label{sg:duality}
\end{align}
Using \eqref{sg:gradientinverse} and $\norm z_2\le N\norm h_2$, we obtain
\[
 |\ip v h|\le\norm{\nabla f}_2\norm{\nabla z}_2+|s|\norm g_2\norm z_2
 \le C K_s\bigl(\norm{\nabla f}_2+\norm g_2\bigr)\norm h_2.
\]
Taking the supremum over unit vectors $h$ and setting $F=\norm{\nabla f}_2+\norm g_2$ gives
\begin{equation}\label{sg:vbound}
 \norm v_2\le C K_sF.
\end{equation}
The identity $isu=v+f$ and Poincar\'e's inequality imply
\begin{equation}\label{sg:ubound}
 |s|\norm u_2\le\norm v_2+\norm f_2\le C K_sF.
\end{equation}
Set $h_0=g+(is+W)f$, so $\norm{h_0}_2\le C_W|s|F$. Taking the real part of the inner product of $P_W(s)u=h_0$ with $u$ and using \eqref{sg:ubound}, we get
\[
 \norm{\nabla u}_2^2\le s^2\norm u_2^2+\norm{h_0}_2\norm u_2
 \le C_W(K_s^2+K_s)F^2\le C_WK_s^2F^2.
\]
Combining this with \eqref{sg:vbound} and $F\le\sqrt2\norm{([f],g)}_{\cH_W}$ proves \eqref{sg:comparebound}. Substitution of \eqref{sg:scalarinput} gives \eqref{sg:generatorlittleo}.
\end{proof}

\subsection{A majorant of positive increase}
Replacing \eqref{sg:generatorlittleo} by $O(s^2)$ and applying the polynomial semigroup theorem \cite{BT} gives only $O(t^{-1/2})$. To retain the strict improvement, we use the general decay theorem of Rozendaal--Seifert--Stahn \cite{RSS}. A continuous nondecreasing function $M:[0,\infty)\to(0,\infty)$ has positive increase if there exist $\alpha,c,S_0>0$ such that
\begin{equation}\label{sg:positiveincrease}
 M(\lambda S)\ge c\lambda^\alpha M(S),\qquad \lambda\ge1,\quad S\ge S_0.
\end{equation}

\begin{lemma}\label{sg:envelope}
If a continuous function $q:\R\to(0,\infty)$ satisfies $q(s)=o(s^2)$, then there is a continuous positive nondecreasing function $M$ of positive increase such that
\begin{equation}\label{sg:envelopemain}
 q(s)\le M(|s|),\qquad M(S)=o(S^2).
\end{equation}
\end{lemma}
\begin{proof}
First set
\[
 g(S)=\max\left\{1,\max_{|r|\le S}q(r)\right\},\qquad S\ge0.
\]
Continuity of $q$ makes $g$ continuous and nondecreasing, with $g(|s|)\ge q(s)$. Given $\varepsilon>0$, choose $R$ such that $q(r)\le\varepsilon r^2$ for $|r|\ge R$, and put $C_R=\max\{1,\max_{|r|\le R}q(r)\}$. For $S\ge R$, we have $g(S)\le C_R+\varepsilon S^2$, so $g(S)=o(S^2)$.

For $S\ge1$, define
\begin{equation}\label{sg:envelopeconstruction}
 H(S)=\max_{1\le r\le S}\frac{g(r)}r,\qquad M(S)=SH(S),
\end{equation}
and set $M(S)=M(1)$ on $[0,1]$. Then $M$ is continuous, positive, and nondecreasing, and $M(S)\ge g(S)$. For $S\ge1$ and $\lambda\ge1$,
\[
 \frac{M(\lambda S)}{M(S)}=\lambda\frac{H(\lambda S)}{H(S)}\ge\lambda,
\]
which proves positive increase. Finally, given $\varepsilon>0$, choose $R\ge1$ such that $g(r)\le\varepsilon r^2$ for $r\ge R$, and set $K_R=\max_{1\le r\le R}g(r)/r$. Splitting the maximum over $[1,R]$ and $[R,S]$ gives
\[
 H(S)\le\max\{K_R,\varepsilon S\},\qquad
 \frac{M(S)}{S^2}\le\max\{K_R/S,\varepsilon\},\qquad S\ge R.
\]
Letting $S\to\infty$ and then $\varepsilon\downarrow0$ proves $M(S)=o(S^2)$.
\end{proof}

Positive increase makes $M$ unbounded. For $t\ge M(0)$, define the generalized inverse
\begin{equation}\label{sg:inverse}
 M^{\leftarrow}(t)=\sup\{S\ge0:M(S)\le t\}.
\end{equation}
Continuity, monotonicity, and unboundedness ensure that the supremum is finite and
\begin{equation}\label{sg:inverseidentity}
 M(M^{\leftarrow}(t))=t.
\end{equation}
Indeed, the set $\{S:M(S)\le t\}$ is closed and bounded. If the value of $M$ at its right endpoint were less than $t$, continuity would allow an extension to the right, contradicting the definition of the supremum.

\subsection{Semigroup decay and completion of the main theorem}
By \cite[Theorem 3.2]{RSS}, if $A$ generates a bounded $C_0$ semigroup $T(t)$ on a Hilbert space, $i\R\subset\rho(A)$, and a continuous nondecreasing function $M$ of positive increase satisfies $\norm{(is-A)^{-1}}\le M(|s|)$, then
\begin{equation}\label{sg:rss}
 \norm{T(t)A^{-1}}=O\left(\frac1{M^{\leftarrow}(t)}\right),\qquad t\to\infty.
\end{equation}

\begin{theorem}\label{sg:decay}
For every fixed $W$ satisfying \eqref{eq:dampingassumptions},
\begin{equation}\label{sg:littleodecay}
 \lim_{t\to\infty}t^{1/2}\norm{T_W(t)G_W^{-1}}_{\cH_W\to\cH_W}=0.
\end{equation}
Consequently, the energy estimate \eqref{eq:mainenergy} in Theorem~\ref{thm:main} holds.
\end{theorem}
\begin{proof}
Lemma~\ref{sg:maximal} and Proposition~\ref{sg:comparison} show that $T_W(t)$ is a contraction semigroup and $i\R\subset\rho(G_W)$. The resolvent identity gives local operator-norm continuity of $s\mapsto(is-G_W)^{-1}$, hence continuity of $q_W(s)$. Applying Lemma~\ref{sg:envelope} to \eqref{sg:generatorlittleo} produces
\[
 \norm{(is-G_W)^{-1}}\le M(|s|),\qquad M(S)=o(S^2),
\]
with $M$ satisfying all hypotheses of \eqref{sg:rss}. Therefore, for sufficiently large $t$,
\begin{equation}\label{sg:inverserate}
 \norm{T_W(t)G_W^{-1}}\le C_W/M^{\leftarrow}(t).
\end{equation}
Set $y=M^{\leftarrow}(t)$. As $t\to\infty$, we have $y\to\infty$ and, by \eqref{sg:inverseidentity}, $t=M(y)$. Hence
\[
 \frac{1/M^{\leftarrow}(t)}{t^{-1/2}}
 =\frac{\sqrt t}{M^{\leftarrow}(t)}
 =\left(\frac{M(y)}{y^2}\right)^{1/2}\longrightarrow0,
\]
which proves \eqref{sg:littleodecay}.

For $U_0=([u_0],u_1)\in D(G_W)$, the state $U(t)=T_W(t)U_0=([u(t)],v(t))$ is differentiable. Choose the representative
\[
 u(t)=u_0+\int_0^t v(\tau)\dd\tau.
\]
Then $u_t=v$ and $v_t=\Delta u-Wv$, recovering the original wave equation with the prescribed initial data. Since $U_0=G_W^{-1}G_WU_0$,
\[
 E(u,t)^{1/2}=2^{-1/2}\norm{T_W(t)U_0}_{\cH_W}
 \le2^{-1/2}\norm{T_W(t)G_W^{-1}}\norm{G_WU_0}_{\cH_W}.
\]
Moreover,
\begin{align*}
 \norm{G_WU_0}_{\cH_W}
 &\le\norm{\nabla u_1}_2+\norm{\Delta u_0-Wu_1}_2\\
 &\le C_W\bigl(\norm{[u_0]}_{\dot H^2}+\norm{u_1}_{H^1}\bigr).
\end{align*}
Thus we may take $\varepsilon_W(t)=2^{-1/2}C_Wt^{1/2}\norm{T_W(t)G_W^{-1}}$. This function tends to zero and gives \eqref{eq:mainenergy} for all $t\ge1$.
\end{proof}

Together with Section~\ref{sec:step6}, this completes the proof of the main theorem. The little-$o$ estimate is uniform in the initial data for fixed $W$; the energy itself satisfies the corresponding $o(t^{-1})$ estimate.

\end{document}